\newif\ifsmfart
\IfFileExists{smfart.cls}
   {\documentclass[12pt,english]{smfart}
    \IfFileExists{smfenum.sty}{\usepackage{smfenum}}{}
    \usepackage{bull}
    \smfarttrue}
{\message{^^J*** It would be better to typeset this file with smfart.cls ***^^J^^J}
\documentclass[12pt]{amsart}}

\usepackage{amssymb}
\usepackage{epsfig}
\usepackage{graphicx}
\numberwithin{equation}{section}

\input xy
\xyoption{all}

\advance\headheight 2pt
\calclayout
\allowdisplaybreaks[3]

\theoremstyle{plain}
\newtheorem{prop}[subsection]{Proposition}

\newtheorem{theo}[subsection]{Theorem}
\newtheorem{coro}[subsection]{Corollary}

\newtheorem{lemm}[subsection]{Lemma}
\newtheorem{assu}[subsection]{Assumption}
\newtheorem{defi}[subsection]{Definition}

\theoremstyle{definition}

\theoremstyle{remark}
\newtheorem{rema}[subsection]{Remark}

\newtheorem{exam}[subsection]{Example}

\newtheorem{nota}[subsection]{Notations}

\def\cJ{{\mathcal J}}

\def\lra{{\longrightarrow}}

\def\Ker{{\rm Ker}}

\newcommand{\La}{\Lambda}
\newcommand{\la}{\lambda}
\newcommand{\al}{\alpha}
\newcommand{\C}{\Bbb C}
\newcommand{\Q}{\Bbb Q}
\newcommand{\Z}{\Bbb Z}

\def\no{\noindent}

\def\cE{{\mathcal E}}
\def\cF{{\mathcal F}}
\def\cG{{\mathcal G}}

\def\cT{{\mathcal T}}

\def\SL{{\rm SL}}
\def\Alb{{\rm Alb}}
\def\ra{\rightarrow}

\def\C{{\mathbb C}}

\def\P{{\mathbb P}}
\def\Q{{\mathbb Q}}

\def\Z{{\mathbb Z}}
\def\C{{\mathbb C}}
\def\N{{\mathbb N}}

\def\codim{{\rm codim}}

\def\Pic{{\rm Pic}}

\makeatother
\makeatletter

\def\Title    { }
\def\Author   {Fedor Bogomolov and Yuri Tschinkel}
\def\Subject  {Algebraic geometry}
\def\Keywords {elliptic fibration , logarithmic transformation}

\newif\ifpdf
   \ifx\pdfoutput\undefined
   \pdffalse
   \else
           \pdfoutput=1
            \pdftrue
   \fi
      \ifpdf
\usepackage[pdftex,
hyperindex=true,
pdftitle={\Title},pdfauthor={\Author},pdfsubject={\Subject},
pdfkeywords={\Keywords},pdfpagemode={UseOutlines},
bookmarks=true,bookmarksopen=true,
bookmarksopenlevel=0,bookmarksnumbered=true,
pdfstartview={FitV},
colorlinks=true
]{hyperref}
 \else
 \usepackage{hyperref}
 \fi

\usepackage{graphicx}

\author{Fedor Bogomolov}
\address{Courant Institute of Mathematical Sciences, N.Y.U. \\
 251 Mercer str. \\
 New York, NY 10012, U.S.A.}
\email{bogomolo@cims.nyu.edu}

\author{Yuri Tschinkel}
\address{Department of Mathematics \\
         Princeton University\\
         Fine Hall, Washington Road\\
         Princeton, NJ 08544-1000,  U.S.A.}
\email{ytschink@math.princeton.edu}

\title[Special elliptic fibrations]
{Special elliptic fibrations}

\begin{document}

\date{\today}

\begin{abstract}
We construct examples of elliptic fibrations
of orbifold general type (in the sense of Campana)
which have no \'etale covers dominating a variety of general type.
\end{abstract}

\maketitle
\tableofcontents

\centerline{
\em To the memory of our friend and colleague Andrey Tyurin.}

\setcounter{section}{0}
\section{Introduction}
\label{sect:introduction}

Consider the following two classes of varieties:
\begin{itemize}
\item admitting an \'etale cover which dominates a
(positive dimensional) variety of general type;
\item admitting a nonconstant map
with target an {\em orbifold} of general type
(defined by taking into account
possible multiple fibers of the map,
see Section~\ref{sect:gen} for details).
\end{itemize}
In this note we construct examples of complex three-dimensional varieties
in the second class which are not in the first class, answering a question
of Campana (see \cite{campana}).

\

\noindent
{\bf Acknowledgments.}
We have benefited from conversations with F. Campana and J. Koll\'ar.
Both authors were partially supported by the NSF.

\

\section{Generalities}
\label{sect:gen}

Throughout, let $X$ be a smooth projective algebraic variety over $\C$
with function field $\C(X)$,
$\Pic(X)$ its Picard group and $K_X$ its canonical class.
For $D\in\Pic(X)_{\Q}$ we let
$\kappa(D)$ be the Kodaira dimension of $D$,
$\kappa(X)=\kappa(K_X)$
the Kodaira dimension of $X$ and $\kappa(X,D):=\kappa(K_X+D)$ the
corresponding log-Kodaira dimension.
We denote by $\Omega^n_X$ the
sheaf of differential  $n$-forms, by $\cT_X$ the tangent bundle
and by $\pi_1(X)$ the fundamental group.

We now recall some notions concerning fibrations following \cite{campana}.
Let $\varphi\,:\,X\ra B$ be a morphism
between smooth algebraic varieties, such that
the locus $D:=\cup_j D_j\subset B$ over which the scheme-theoretic fibers
of $\varphi$
are not smooth is a (strict) normal crossing divisor
(with irreducible components $D_j$).
For each $j$, let $n_j$ be the minimal (scheme-theoretic) multiplicity of
a fiber-component over $D_j$ and
$$
D(\varphi):=\sum_j (1-1/n_j)D_j\in \Pic(B)_{\Q}
$$
the multiplicity divisor of $\varphi$.
The pair $(B,D(\varphi))$ will be called an {\em orbifold}
associated to $\varphi$.
It is called of an orbifold of {\em general type} if
$$
\kappa(B,D(\varphi))=\dim(B)>0.
$$

\begin{exam}
\label{exam:dim2}
Let $\varphi\,:\, X\ra B=\P^1$ be an elliptic fibration
such that $D(\varphi)\neq \emptyset$.
The degenerate fibers with $n_j\ge 2$ are {\em multiple}
fibers. The associated orbifold $(B,D(\varphi))$ is
of general type provided
\begin{equation}
\label{eqn:gene}
\sum (1-1/n_j)>2.
\end{equation}
This condition implies that
there exists a finite cover $\tilde{B}\ra B$ ramified
with multiplicity $n_j$ at points $D_j\subset D(\varphi)$
and of genus $\ge 2$.
Let $\tilde{X}$ be the pullback of the elliptic fibration to $\tilde{B}$.
Then $\tilde{X}\ra X$ is \'etale and 
has a surjective map $\tilde{X}\ra\tilde{B}$.
\end{exam}

\begin{theo}
\label{thm:main}
There exist smooth projective algebraic threefolds $X$ admitting 
an elliptic fibration 
$\varphi\,:\, X\ra B$ such that 
\begin{itemize}
\item $\pi_1(X)=0$;
\item $B$ is a smooth elliptic surface with $\kappa(B)=1$;
\item $D(\varphi)\subset B$ is a smooth irreducible divisor;
\item the orbifold $(B,D(\varphi))$ is of general type.
\end{itemize}
\end{theo}

\section{Logarithmic transforms}
\label{sect:log}

We recall the construction of
logarithmic transforms of elliptic fibrations due to Kodaira
\cite{kodaira}
(for more details see \cite{FM}, Section 1.6).

Let $C$ be a smooth curve and
$\eta\,:\,\cE\ra C$ a {\em nonisotrivial} elliptic fibration.
Let $\Delta \subset C$ be a unit disc with center
$p_0$ and {\em smooth} central fiber $E_0$ over $p_0$.
For every $m\in \mathbb N$ consider the diagram

\centerline{
\xymatrix{
   \tilde{\cJ} \ar[d]_{\tilde{\eta}} \ar[r] &  \ar[d]^{\eta} \cJ \\
   \tilde{\Delta}     \ar[r]_{\iota_m}      &   \Delta
         }
}

\no
where  $\cJ$ is
the restriction of $\cE$ to $\Delta$,
$\iota_m$ is a cyclic cover of degree $m$ given by
$$
\tilde{z}\mapsto \tilde{z}^m=z
$$
(with $z$ a local analytic coordinate at $p_0$) and
$\tilde{\cJ}$ the pullback of $\cJ$ to $\tilde{\Delta}$.
After appropriate choices one has
$$
\cJ=(\C\times \Delta)  /\La(z),  \,\,\,\,\,\,\, \tilde{\cJ}=(\C\times
\Delta)  /\La(\tilde{z}^m)
$$
(where $\La(z)\subset \C$ is a family of lattices)
and
$$
\cJ=\tilde{\cJ}/\mathfrak C_m,
$$
where $\mathfrak C_m$ is a finite cyclic group generated by
$$
(s,\tilde{z})\mapsto (s,\zeta_m \tilde{z}) \,\,\mod \La(\tilde{z}^m)
$$
(and $\zeta_m$ is an $m$-th root of 1).
Let $\omega_m(z)/m$ be a local $m$-torsion section of $\cJ$ and
define
$$
\cJ':= \tilde{\cJ}/{\mathfrak C}'_m,
$$
where ${\mathfrak C}'_m$ is a cyclic group generated by
$$
(s,\tilde{z})\mapsto
(s +\frac{\omega_m(\tilde{z}^m)}{m},\zeta_m \tilde{z}) \,\,\mod
\La(\tilde{z}^m).
$$
We have an isomorphism
$$
\cJ'\setminus (E_0 /\mathfrak C_m) \simeq \cJ\setminus E_0
$$
and we can extend $\cJ'$ to an elliptic fibration 
$\eta'\,:\, \cE'\ra C$, 
called the {\em logarithmic transform} (twist) of $\cE$. 
In $\cE'$ a cycle (circle) $\mathbb S$ which
was bounding a holomorphic section
over a disc in $\cE$ is homologous to a nontrivial
cycle $\mathbb S'\in E_0$.

\begin{prop}
\label{prop:alg}
Assume that $\cE$ is locally Jacobian and not locally
isotrivial and that $\cE'$ is obtained from $\cE$ by a 
logarithmic transform at exactly one point $p_0\in C$.
Then 
\begin{itemize}
\item  
$H^0(\cE,K_{\cE})\simeq H^0(\cE',K_{\cE'})$;
\item 
$\pi_1(\cE ) = 0 \Rightarrow \pi_1(\cE')=0$;
\item 
$\cE'$ is K\"ahler.
\end{itemize}
\end{prop}

\begin{proof}
Every form $w\in H^{2,0}(\cE)$ has a local representation as
$$
w=dh\wedge d \log(s).
$$ 
It is visibly invariant under translation by $s$ on $\cE\setminus E_0$, is
preserved under gluing and can be extended from 
$\cE\setminus E_0$ to $\cE'$. 
Moreover, on $\cE'$ it has a 
zero of multiplicity $m-1$ along $E_0/\mathfrak C_m$.
After twisting exactly one fiber,
we have
$$
K_{\cE'} = K_{\cE} + (1-1/m) E,
$$
where $E$ is a (generic) fiber of $\eta$.
Since $\cE$ is locally Jacobian we have  $K_{\cE} = \eta^* L$,
where $L\in \Pic(C)$, and $H^0(\cE, K_{\cE}) = H^0(C,L)$. Similarly, we have 
an imbedding $K_{\cE'} \hookrightarrow \eta^*(L+ p_0)$ and 
$$
H^0(\cE',K_{\cE'})\subset H^0(\cE', {\eta'}^*(L+ p_0)) = H^0(C,(L+ p_0)).
$$
We have 
$h^0(C,(L+ p_0)) \leq h^0(C,L) + 1$. An 
$f\in H^0(C,(L+ p_0))$ which is not in the image 
of $H^0(C,L)$ is nonzero at $p_0$. The corresponding element
in $H^0(\cE', {\eta'}^*(L+ p_0))$ is also nonzero on the fiber over $p_0$.
However, every global section of  $K_{\cE'}$ vanishes on the 
central fiber. Thus $f\notin H^0(\cE',K_{\cE'})$ 
so that every section of $K_{\cE'}$ 
is an extension of a section of $K_{\cE}$ 
(restricted to $\cE\setminus E_0$): 
$$
H^0(\cE, K_{\cE}) =H^0(\cE',K_{\cE'}).
$$
Since $\pi_1(\cE) = 0$ we have $C=\P^1$.
We claim that $\pi_1(\cE \setminus E_0) = 0$. Indeed,
the fundamental group of the elliptic fibration 
$(\cE \setminus E_0)\ra (C\setminus p_0)$ 
lies in the image of $\pi_1(E)$, which is 
a {\em finite} abelian group since there are nontrivial vanishing cycles
which are homotopic to zero. Since the global monodromy has 
finite index in $\SL_2(\Z)$ the 
group $\pi_1(\cE \setminus E_0)$ is also finite abelian
and the corresponding covering is fiberwise. 
Thus it extends as a finite \'etale covering of $\cE$, contradicting 
the assumption that $\pi_1(\cE) = 0$.

Consider the (topological) quotient spaces $\cE/ E_0$ and $\cE'/E'_0$.
They are naturally isomorphic and we have two exact homology sequences
$$
\ldots 
\ra H_3(\cE, \Q)\ra H_3( \cE/E_0,\Q)
\stackrel{d}{\lra} H_2(E_0,\Q) \ra H_2(\cE, \Q)
$$
and 
$$
\ldots 
\ra H_3(\cE', \Q)\ra H_3( \cE'/E_0',\Q)
\stackrel{d'}{\lra} H_2(E_0',\Q) \ra H_2(\cE',\Q).  
$$   
Since $\cE$ is K\"ahler
$$
\Q=  H_2(E_0,\Q)\hookrightarrow H_2(\cE,\Q)
$$
and 
$$
H_3(\cE, \Q) = H_3(\cE/E_{0},\Q) = H_1(C,\Q)^*.
$$
Here we used that $H_1(\cE,\Q) = H_1(C,\Q)$ which follows from the
local nonisotriviality of $\cE$. 
Geometrically it means that every $3$-cycle on $\cE$ and $\cE/E_0$
can be realized as a product of a $1$-cycle on $C$ and an elliptic fiber.

Since $\cE/E_{0} = \cE'/E'_{0}$ the differential $d'$ is also zero
and 
$$
H_2(E_{0}',\Q) \hookrightarrow  H_2(\cE',\Q).
$$
Thus the class of the generic fiber $E$ is nontrivial.
This implies the existence of a K\"ahler metric (see \cite{miyaoka}).
Therefore, if $\cE$ is K\"ahler then so is $\cE'$.
\end{proof}

\begin{coro}
\label{coro:alg}
If $\cE$ is algebraic and rational then $\cE'$ is algebraic.
\end{coro}

\begin{proof}
A smooth surface $S$ is projective iff 
there is a class $x\in H_2 (S,\Q)$ with $x^2 > 0$ which 
is orthogonal to $H^{2,0}(S)\subset H^2(S,\C)$.
Since $\cE'$ is K\"ahler and $H^{2,0}(\cE) = H^{2,0}(\cE') = 0$
there is such a class in $H_2(\cE, \Q)$.
\end{proof}

\begin{exam}
\label{exam:hi}
Let $\xi\, :\,\P^1 \ra \P^1$ be a polynomial map of degree $n\ge 2$
which is cyclically $n$-ramified over $\infty$.
Let $\bar{\varphi}\,:\, \bar{\cE}\ra \P^1$ 
be a rational elliptic surface and $\bar{\cE}'$ 
its logarithmic $nm$-twist over $\infty$. 
Consider the diagram

\centerline{
\xymatrix{
\cE\ar[d]_{\xi}\ar[r]^{\eta} & \P^1  \ar[d]_{\xi} & \ar[l]_{\eta'} \ar[d]^{\xi}\cE'\\     
\bar{\cE} \ar[r]_{\bar{\eta}} & \P^1  & \ar[l]^{\bar{\eta}'} \bar{\cE}'
          }
}

\

\noindent
The surface $\cE'$ (induced by $\xi$) is a logarithmic $m$-twist at $\infty$
of $\cE$ (induced from $\bar{\cE}$). We have
$h^{0} (\cE,K_{\cE}) = n - 1$. Since $\bar{\cE}'$ is algebraic 
(by Corollary~\ref{coro:alg}), $\cE'$ is also algebraic.     
\end{exam}

For more details concerning algebraicity of elliptic fibrations
obtained by logarithmic transformations we refer to
\cite{FM}, Section 1.6.2.

\section{Construction}
\label{sect:const}

Consider the following diagram

\centerline{
\xymatrix{
X  \ar[d]_{\varphi} &                                       
&  \ar[d]^{\eta'}\cE'    \\
B  \ar[r]_{\beta}   &\ar[d]^{\varphi_1}\ar[r]_{\psi} S         & C                    \\
                    & C_1                                   &
          }
}

\

\no
where
\begin{itemize}
\item $C_1,C$ are $\P^1$;
\item $S$ is a nonisotrivial locally Jacobian elliptic
surface with  irreducible fibers, $\pi_1(S)=0$ and $\kappa(S) = 1$;
\item $\psi \,:\, S \to C$ is a rational map
with connected fibers defined by a generic
line $\P^1_\psi\subset \P(H^0(S,L))$, where $L$ is a polarization on $S$;
\item $\beta\, :\, B\ra S$ is a minimal blowup 
so that $\gamma:=\psi\circ\beta\, :\, B\ra C$ 
is a fibration with irreducible
fibers (it exists since $L$ is very ample
and $\psi$ is generic, i.e., all singularities of $\psi$ are simple
and lie in different smooth fibers of $\varphi_1$);
\item $\eta'\,:\, \cE'\ra C$ is the fibration from Example~\ref{exam:hi};
\item $\varphi\,:\, X\ra B$ is the pullback of $\eta'$ via $\gamma$.
\end{itemize}

\begin{lemm}
\label{lemm:gg}
Let  $B$ be the surface above, $p\in C$ a generic point and
$D = \gamma^{-1}(p)$.
Then
\begin{itemize}
\item $\varphi_1\circ \beta\,:\, B\ra C_1$ 
is an elliptic fibration and $\kappa(B) = 1$;
\item $\pi_1(B \setminus D) = 0$.
\end{itemize}
\end{lemm}

\begin{proof}
The genericity of $L$ and $\psi$ implies that all  
fibers are irreducible and that $B$ is
a blowup of $S$ in a finite number of distinct points 
in which the divisors from $\P^1_\psi\subset \P(H^0(S,L))$
intersect transversally.
Since $\pi_1(S)=0$ we have
$$
\pi_1( S \setminus D) = \mathfrak C_m
$$ 
(by Lefschetz theorem), where $m$ is the largest
integer dividing $L$ in $\Pic(S)$.
The corresponding cyclic
covering of $S$ is $m$-ramified along $D$.
We have
$$
(B \setminus D)=(S\setminus D) \cup \bigcup_{i\in I} \ell_i
$$ 
where $I$ is a finite set and $\ell_i$ are affine lines.
A cycle generating 
$\pi_1(B\setminus D)$ is contracted inside one of these lines,
which implies that the image of $\pi_1(S\setminus D)$ in 
$\pi_1(B\setminus D)$ is trivial.
\end{proof}

\

\begin{proof}[Proof of Theorem~\ref{thm:main}]
The elliptic fibration $\varphi\,:\, X\ra B$
satisfies the claimed properties.

First observe that $D$ intersects all components of
the fibers $E$ of  
$$
\varphi_1\circ \beta\,:\, B\ra C_1.
$$ 
Indeed, by genericity every such $E$ 
is either irreducible or a union
of a smooth elliptic curve and a rational $(-1)$ curve $P$.
If $E$ is irreducible the claim follows
from the ampleness of $L$. For the same reason we have 
$\deg(D_{|E}) \geq 2$. Since $D\cdot P = 1$ 
there is a nontrivial intersection
with another component.

Put $F:=K_B + (1-1/m)D$. 
Since $\kappa(B) =1$ and $K_B\cdot D > 0$
a subspace of sections in $H^0(B,amF)$ (for some $a\in \N$) 
gives a surjection $B\ra C_1\times C$, so that $\kappa (F) = 2$.
Moreover, $F$ intersects positively every divisor in
$B$ (except finitely many rational curves $P_i$ obtained
by blowing up $S$).
It follows that $F=H+\sum m_iP_i$, 
where $H$ is a polarization on $B$ and $m_i\ge 0$. 
Thus $(B,D(\varphi))$ is an orbifold of general type. 
In particular, 
$$
\varphi^* K_B\subset \varphi^*F\subset\Omega^2_X,
$$
where $\varphi^*F$ is saturated, and $\kappa(\varphi^*F)=2$.
Notice that $\kappa(X)=2$ since
$$
\varphi^*F\times \gamma^*K_{\cE'/C}\subset K_X\,\,\, 
\text{ and } \,\,\, \kappa(K_{\cE'/C}) = 1.
$$
The pullback $H'$ (to $X$) of a polarization on $\cE'$ 
is positive on the fibers of $\varphi$.
Since $B$ is projective it has a polarization $H$ such that
for some $a\in \N$ the divisor 
$a\varphi^*H + H'$ is positive on every curve in $X$ 
and is represented by a positive definite K\"ahler form. 
This implies that $X$ is projective.

We claim that $\pi_1(X) = 0$. We know that 
$\pi_1(B) =\pi_1(B\setminus D)=0$. Hence
$\pi_1(X)$ is in the image of $\pi_1(E)$ of a smooth (elliptic)
fiber of $\varphi$. Since the monodromy of $\varphi$ is large it
kills the fundamental group of the fiber. Indeed, the restriction
of $\varphi$ to a $\P^1$ is isomorphic to $\cE'$. 
Since the complement of a multiple fiber
in $\cE'$ has trivial fundamental group the same holds
for $X$.

Thus $X$ admits a map onto an orbifold of general type
but does not dominate a variety of general type nor has (any) \'etale
covers.
\end{proof}

\begin{rema}
In fact, we have proved that  $\pi_1(X \setminus \varphi^{-1 }(D)) = 0$
so that no modification can yield an \'etale cover.
\end{rema}

\section{Holomorphic differentials}
\label{sect:sub}

One of the features of the construction in Section~\ref{sect:const}
was the use of a 1-dimensional subsheaf of holomorphic
forms with many sections. We have seen that such sheaves 
impose strong restrictions on the global geometry of the variety. 
Generalizing several results in \cite{bo}, we now
give an alternative proof of Campana's theorem 
on the correspondence between such sheaves and maps
onto orbifolds of general type (see \cite{campana}).

Let $X$ be a smooth K\"ahler manifold and $\omega\in \Omega_X^i$
a form. The {\em kernel} of $\omega$ is the subsheaf of
$\cT_X$ generated (locally) by sections $t$  
such that for all $x\in \Lambda^{i-1}\cT_X$
$$
\omega(t\wedge x)=0.
$$
The kernel doesn't change under multiplication of $\omega$ by 
a nonzero (local) holomorphic section of the structure sheaf.
This defines, for every subsheaf  $\cF\subset \Omega^i_X$, 
its kernel  $\Ker({\cF})$ (a special case of the notion 
of support of a differential ideal).

\begin{defi} 
\label{defi:su}
We say that $\cF\subset \Omega^i_X$ 
is $k$-monomial if 
at the generic point of $X$ a
nonzero local section $f$ of $\cF$ 
is a product of local holomorphic $1$-forms:
$$
f=q_1\wedge \ldots \wedge q_k \wedge \omega,
$$ 
where $1\le k \le i$ and $\omega$ is a local $(i-k)$-form.
We call $\cF$ monomial if $k=i$.  
\end{defi}

\begin{prop} 
\label{prop:1} 
Let $X$ be a smooth compact K\"ahler manifold
and $\cF\subset \Omega^{i}_X$ a one-dimensional
subsheaf such that 
$$
h^0(X, \cF^n ) \geq an^{k} + b,
$$ 
where $a >0 $ and $k\ge 1$. Then 
\begin{itemize}
\item $k\le i$;
\item $\cF\subset \Omega_X^i$ 
is a $k$-monomial subsheaf;
\item  there exist an algebraic variety $Y$ of 
dimension $k$ and a meromorphic map 
$$
\varphi=\varphi_{\cF} \,:\, X \ra Y
$$
with irreducible generic fibers 
such that the tangent space of the fiber of 
$\varphi$ at a generic point coincides 
with $\Ker(\cF)$.
\end{itemize}
\end{prop}

\begin{proof}
The ratios of sections $s_l\in H^0(X,\cF^n)$ 
generate a field of transcendence degree $k$ (for some $n\ge 1$).
In particular, there is an $x\in X$, with $s_0(x)\neq 0$, 
where the local coordinates
$$
f_l = s_l(x)/s_0(x), \,\,\,  l=1,...,k
$$
are independent. We know that $s_0$ is locally equal to $w_0^n$,
where $w_0$ is a local closed form nonvanishing at $x$  (see \cite{bo}).
Further,  $f_l w_0$ is also a local closed form nonvanishing at $x$. Since 
$$
ds_0 = d(f_l w_0) = df_l\wedge w_0 = 0
$$
we obtain  
$$
w_0 = df_l\wedge w'.
$$ 
Since the forms $df_l$ are linearly independent we see that
\begin{itemize}
\item  $w_0 = g df_1\wedge df_2 \ldots \wedge df_k\wedge\omega $, so that $g$ is
algebraically dependent on $f_l$ and $\cF$ is a  $k$-monomial subsheaf of $\Omega^i_X$;
\item the fibers of the map given by $f_l$ are
tangent to the kernel of $w_0$.
\end{itemize}
Thus we have a meromorphic map 
\begin{equation}
\label{eqn:var}
\varphi\, :\, X\ra Y,\,\,\, \dim(Y)=k,
\end{equation}
such that $s_l$ are locally (at a generic point of $X$) products
of a power of a volume form induced from $Y$  under
$\varphi$ and a power of $\omega$ which
is nontrivial on the fiber of $\varphi$.
The map $\varphi$ is holomorphic outside of the zero locus
of the ring $\oplus_n H^0(X,\cF^n)$.
\end{proof}

\begin{coro}
\label{coro:mono}
If $k = i$ or $k=i -1$ then $\cF\subset\Omega^i_X$ is monomial.
\end{coro}

\begin{proof}
It suffices to consider $f\in \cF$ at
generic points. There are two cases:
\begin{itemize}
\item $k =i$: then
$$
f= df_1\wedge \ldots \wedge df_k,
$$ 
(modulo multiplication by a function).
\item  
$k = i-1$: then   
$$
f = df_1\wedge df_2 .......\wedge df_k\wedge q,
$$
where $q$ is a closed $1$-form.
\end{itemize}
\end{proof}

\begin{rema}
\label{rem:map}
The map from \eqref{eqn:var} admits a bimeromorphic modification
$$
\varphi\,:\,X\ra Y
$$ 
such that
\begin{itemize}
\item $\varphi$ is holomorphic with generically smooth
and irreducible fibers;
\item $X$ and $Y$ are smooth.
\end{itemize}
\end{rema}

\begin{nota}
\label{nota:dege}
For $\varphi$ as in Remark~\ref{rem:map}
we define its {\em degeneracy locus} $D=D_\varphi$ 
as the subset of all $y\in Y$ such that 
$d\varphi(x)=0$ for all $x\in \varphi^{-1}(y)$.  
\end{nota}

\begin{rema}
\label{rem:furt}
After another modification of $\varphi$ we can achieve that 
\begin{itemize}
\item $\codim(D)\ge 2$ or
\item $D=\cup_j D_j$ and 
each $\tilde{D}_j:=\varphi^{-1}(D_j)=\cup_i \tilde{D}_{ij}$ 
is a normal crossing divisor.
\end{itemize}
\end{rema}

\begin{assu}
\label{assu:gen}
The map $\varphi$ is as in Remarks~\ref{rem:map} and \ref{rem:furt}.
\end{assu}

\begin{lemm} 
\label{lemm:orb}
If $k=i$ then either $\codim(D)\ge 2$ and 
$$
\cF=\varphi^* K_Y
$$
or there exist integers $n_j\ge 1$ such that 
\begin{itemize}
\item $D(\varphi):=K_Y + \sum_j (1 - 1/n_j) D_j$ is big on $Y$;
\item $\varphi$ has multiplicity $\ge n_j$ along every
$\tilde{D}_{ij}$;
\item $\varphi^* D(\varphi)\subset \cF$.
\end{itemize}
\end{lemm}

\begin{proof}
Every $x$ with $d\varphi(x)\neq 0$
has a neighborhood $U$ such that the restriction of 
every section $s\in H^0(X,\cF^n)$
to $U$ is induced from a (unique) section 
$s_U\in H^0(\varphi(U),nK_Y)$.
There is a unique holomorphic tensor 
$s_Y\in H^0(Y\setminus D,nK_Y)$  
(where $D=D_{\varphi}$ is the degeneracy locus of $\varphi$) 
such that the restriction of 
$\varphi^*(s_Y)$ to $X\setminus \varphi^{-1}(D)$
coincides with $s$.

If $\codim(D) \geq 2$ then $s_Y$
has a unique extension to a holomorphic tensor on $Y$
(since $Y$ is smooth). In this case, $Y$ is of general type. 
In case $\codim(D)=1$  we see (using Remark~\ref{rem:furt})     
that $s_Y$ is a well-defined tensor on $Y$ with poles 
along $D_j$, i.e.,  
$$
s_Y\in H^0(Y,nK_Y+ \sum_j d_j D_j),
$$
for some $d_j\in \N$.
Let $n_j$ be the minimal multiplicity of $\varphi$ on the components 
$\tilde{D}_{ij}$ which surject onto $D_j$ (for all $j$). 
Since $\varphi^* s_Y$ is holomorphic on $X$, a local computation shows that
$$
d_j\le n(1-1/n_j)
$$ 
(see, for example, \cite{sakai} and \cite{bo}) and  that
$$
K_Y+(1-1/n_j)D_j
$$
is big. 
\end{proof}

\begin{rema} 
This gives an alternative proof of 
Campana's theorem characterizing 
fibrations over orbifolds of general type.
\end{rema}

In the case $i =k$ a section of $\cF^n$ (at a generic point of $X$) 
descends to the $n$th-power of a (local) volume form on $Y$ but
the correponding global form on $Y$ may have singularities. 
These singularities disappear after
a finite local covering which is sufficiently ramified along the singular 
locus. This property can be defined for arbitrary tensors.

\begin{defi} 
\label{defi:loca}
A meromorphic tensor
$t$ on $Y$ is {\em locally integrable} if
for every point $y\in Y$ there exist
a neighborhood $U=U_y$ and a (local) manifold $V$ together
with a proper finite map 
$$
\la \,:\,\varphi^{-1}(U)\ra U
$$
such that $\la^* t$ is holomorphic on $V$.
\end{defi}

For $k = i-1$ we have an analog of Lemma~\ref{lemm:orb}:

\begin{lemm} 
\label{lemm:smo}
Let $X$ be a smooth compact K\"ahler manifold,
$\cF\subset \Omega^i_X$ a one-dimensional subsheaf such that
$$
a'n^{i-1} + b' > h^0(X, \cF^n ) \geq an^{i-1}+ b,
$$
with $a >0$, and $\varphi=\varphi_{\cF}$ 
(as in Proposition~\ref{prop:1}).
Then there exist a nontrivial fibration 
$$
\rho \,:\,A_Y \ra Y
$$ 
(with fibers complex tori)
and a map
$$
\al=\al_\cF \,:\, X \ra  A_Y
$$  
with connected fibers such that 
\begin {itemize}
\item  $\varphi = \rho\circ \al$; 
\item  the tangent space of the fiber of $\al$ at a
generic point is contained in $\Ker(\cF)$;
\item there is a divisor $D\subset Y$ such that every section
$s\in \cF $ is a lifting of a monomial locally integrable tensor
on $A_Y$. 
\end{itemize}
\end{lemm}

\begin{proof} 
By Proposition~\ref{prop:1}, there is map $\varphi \,:\, X\ra Y$, where
$\dim(Y)=i-1$.  It has a natural factorization
$$
\varphi \,:\, X\stackrel{\al}{\lra} A\stackrel{\rho}{\lra} Y,
$$
where the fiber $A_y$ of $\rho$ over a generic $y\in Y$ 
is the Albanese variety $\Alb(X_y)$. 
Since $\cF \subset \Omega^i_X$ any section $s\in \cF^n$ 
(at a generic point of $X$) can be represented as 
$$
s=(df_1\wedge..... \wedge q)^n,
$$ 
where $q$ is a closed $1$-form. The form $q$ defines a holomorphic form on
a generic fiber $X_y$ of $\varphi$ so that $\rho$ is {\em nontrivial}
and $\dim(\alpha(X_y))\ge 1$.
In particular, the restriction of $q$ to ${X_y}$ 
is induced from $A_y$.

It follows that there exists a sheaf  $\cG\subset \Omega^i_{A}$
such that $\cF^n$ is a saturation of 
$\al^*\cG$. Moreover, all sections of $\cF^n$
are obtained as lifts of integrable meromorphic sections of $\cG$.
\end{proof}

\begin{rema}
Notice that if $\dim(A_y)=1$ (for generic $y\in Y$)
then $\cG=K_A$ (and $A\ra Y$ is an elliptic fibration). 
\end{rema}

\bibliographystyle{smfplain}
\bibliography{ellhyp}

\providecommand{\bysame}{\leavevmode ---\ }
\providecommand{\og}{``}
\providecommand{\fg}{''}
\providecommand{\smfandname}{\&}
\providecommand{\smfedsname}{\'eds.}
\providecommand{\smfedname}{\'ed.}
\providecommand{\smfmastersthesisname}{M\'emoire}
\providecommand{\smfphdthesisname}{Th\`ese}
\begin{thebibliography}{1}

\bibitem{bo}
{\scshape F.~A. Bogomolov} -- {\og Analytic sections in conic bundles\fg},
  \emph{Trudy Mat. Inst. Steklov.} \textbf{165} (1984), p.~16--23, Algebraic
  geometry and its applications.

\bibitem{campana}
{\scshape F.~Campana} -- {\og Special varieties and classification theory\fg},
  {\tt ArXiv: math.AG/0110051}, 2001.

\bibitem{FM}
{\scshape R.~Friedman {\normalfont \smfandname} J.~W. Morgan} -- \emph{Smooth
  four-manifolds and complex surfaces}, Ergebnisse der Mathematik und ihrer
  Grenzgebiete (3) [Results in Mathematics and Related Areas (3)], vol.~27,
  Springer-Verlag, Berlin, 1994.

\bibitem{kodaira}
{\scshape K.~Kodaira} -- {\og On the structure of compact complex analytic
  surfaces. {I}\fg}, \emph{Amer. J. Math.} \textbf{86} (1964), p.~751--798.

\bibitem{miyaoka}
{\scshape Y.~Miyaoka} -- {\og K\"ahler metrics on elliptic surfaces\fg},
  \emph{Proc. Japan Acad.} \textbf{50} (1974), p.~533--536.

\bibitem{sakai}
{\scshape F.~Sakai} -- {\og Kodaira dimensions of complements of divisors\fg},
  in \emph{Complex analysis and algebraic geometry}, Iwanami Shoten, Tokyo,
  1977, p.~239--257.

\end{thebibliography}

\end{document}